\documentstyle[leqno,12pt,amstex]{article}
\numberwithin{equation}{section}
\begin{document}
\title
{
\large{\bf The global existence theorem for quasi-linear 
wave equations with multiple speeds, II}\\
{~}\\
\small{\it Dedicated to Professor Thomas C. Sideris on 
the occasion of his 60th birthday}
}
\author
{
\large{Kunio Hidano}
}
\date{}
\maketitle
\begin{abstract}
The Cauchy problem is studied for systems of quasi-linear wave equations 
with multiple speeds in two space dimensions. 
Using the method of Klainerman and Sideris together with 
the localized energy estimate, 
we give an alternative proof of a beautiful result of Hoshiga and Kubo.
\end{abstract}


\section{Introduction}
\baselineskip=0.6cm
In his previous paper \cite{Hi4} with the same title as above, 
the present author considered the problem of 
global existence of small solutions 
to the Cauchy problem for systems of quasi-linear wave equations with 
multiple speeds. 
Relying entirely upon the Klainerman-Sideris method \cite{KS}, 
he aimed at giving a unified proof of the beautiful results 
of Hoshiga and Kubo \cite{HK} and Yokoyama \cite{Yo} 
whose proofs built upon point-wise (in space and time) hard estimates 
of the fundamental solutions. 
This aim of \cite{Hi4} was left unaccomplished. 
Indeed, in the case of three space dimensions 
the same result as Yokoyama \cite{Yo} was obtained 
without any reliance upon such hard point-wise estimates of 
the fundamental solution. 
In the case of two space dimensions, however, 
by the Klainerman-Sideris method, 
the author was only able to obtain a global existence result 
under the very restrictive assumption that 
initial data were compactly supported. 
The purpose of revisiting \cite{Hi4} is 
to remove this assumption and 
obtain the same result as Hoshiga and Kubo \cite{HK}. 
We rest our proof on the method of Klainerman and Sideris 
together with two other ingredients: 
the space-time $L^2((0,\infty)\times{\mathbb R}^2)$ estimate and 
the ${\dot H}^{1/2}({\mathbb R}^2)$ estimate. 
In this way, we avoid using the Hardy-type inequality of Lindblad 
(Lemma 1.2 of \cite{Lin}, Lemma 3.3 of \cite{Ka}), 
while in \cite{Hi4} we resorted to using it under the condition that 
the solutions were compactly supported in the space variables 
for every fixed time. 

This paper is organized as follows. We explain the notation 
in the next section, and then state the main theorem in Section 3. 
Useful Sobolev-type inequalities and crucial estimates of the null 
forms are collected in Section 4. 
Weighted $L^2({\mathbb R}^2)$-norms of the second or higher-order derivatives 
are shown to be bounded by generalized energies in Section 5. 
In Section 6, we carry out 
higher-order energy estimates, space-time $L^2$-estimates, 
${\dot H}^{1/2}({\mathbb R}^2)$ estimates, 
and lower-order energy estimates 
to complete the proof of the main theorem. 
\section{Notation}
As in the previous paper \cite{Hi4}, 
we follow Sideris and Tu $\cite{ST}$ to use the notation. 
Repeated indices are summed if lowed and uppered. 
Greek indices range from $0$ to $2$ (space dimensions), 
and roman indices from $1$ to $m$. 
We shall consider 
systems of $m$ quasi-linear equations. 
Points in $(0,\infty)\times{\mathbb R}^2$ are denoted by 
$(x^0,x^1,x^2)=(t,x)$. 
In addition to the usual partial differential operators 
$\partial_\alpha=\partial/\partial x^\alpha$ 
$(\alpha=0,1,2)$ with the abbreviation 
$\partial=(\partial_0,\partial_1,\partial_2)
=(\partial_0,\nabla)$, 
we use the generator of Euclid rotation 
$\Omega=x^1\partial_2-x^2\partial_1$ 
and of space-time scaling $S=x^\alpha\partial_\alpha$. 
The set of these $5$ vector fields is denoted by 
$\Gamma=\{\,\Gamma_0,\Gamma_1,\dots,\Gamma_4\,\}
=\{\,\partial,\Omega,S\,\}$. 
We employ the multi-index notation in Sideris and Tu $\cite{ST}$ 
to mean $\Gamma^a:=\Gamma_{a_\kappa}\cdots\Gamma_{a_1}$ 
for $a=(a_1,\dots,a_\kappa)$, 
a sequence of indices $a_i\in\{0,\dots,4\}$ of length $|a|=\kappa$. 
It is convenient to set 
$\Gamma^a=1$ if $|a|=0$. Suppose that 
$b$ and $c$ are disjoint subsequences of $a$, allowing 
that $|b|=0$ or $|c|=0$. We say 
$b+c=a$ if $|b|+|c|=|a|$, 
$b+c<a$ if $|b|+|c|<|a|$. 

The D'Alembertian, which acts on vector-valued 
functions $u:(0,\infty)\times{\mathbb R}^2\to{\mathbb R}^m$, is denoted by
$$
\square
=
\mbox{Diag}(\square_1,\dots,\square_m),\quad
\square_k=\frac{\partial^2}{\partial t^2}-c_k^2\Delta.
$$
In what follows, 
we suppose that each of the $m$ propagation speeds 
is different from the other $m-1$. 
Without loss of generality, we assume
\begin{equation}
0<c_1<c_2<\cdots<c_m.
\end{equation}

Associated with this operator, 
the standard and generalized energies are defined as 
\begin{eqnarray*}
& &
E_1(u(t))
=
\frac12
\sum_{k=1}^m
\int_{{\mathbb R}^2}
\bigl(
|\partial_t u^k(t,x)|^2
+
c_k^2|\nabla u^k(t,x)|^2
\bigr)
dx,\\
& &
E_\kappa(u(t))
=
\sum_{|a|\leq\kappa-1}
E_1(\Gamma^a u(t)),\quad
\kappa=2,3,\dots
\end{eqnarray*}
Allowing a higher-order energy to grow polynomially in time 
but bounding a lower-order one uniformly in time, we build up 
a series of estimates of the generalized energies. 
The auxiliary norm 
$$
M_\kappa(u(t))
=
\sum_{k=1}^m
\sum_{|a|=2}
\sum_{|b|\leq\kappa-2}
\|\langle c_kt-|x|\rangle
\partial^a\Gamma^b
u^k(t)\|_{L^2({\Bbb R}^2)},\quad \kappa=2,3,\dots
$$
plays an intermediate role. 
Here and later on as well 
we use the notation 
$\langle A\rangle=\sqrt{1+|A|^2}$ for a scalar or a vector 
$A$. 
We also use 
\begin{eqnarray} 
& &
N(u(t))
:=
\biggl(
\sum_{k=1}^m
\bigl(
\||D_x|^{1/2}u^k(t)\|_{L^2({\mathbb R}^2)}
+
\||D_x|^{-1/2}\partial_tu^k(t)\|_{L^2({\mathbb R}^2)}
\bigr)
\biggr)^{1/2},\\
& &
N_\kappa(u(t))
:=
\sum_{|a|\leq \kappa}N(\Gamma^a u(t)).\nonumber
\end{eqnarray}
\section{Result}
We consider the Cauchy problem for a system 
of quasi-linear wave equations
\begin{equation}
\square u=F(\partial u,\partial^2u)
\quad\mbox{in}\,\,(0,\infty)\times{\mathbb R}^2
\end{equation}
subject to the initial data
\begin{equation}
u(0)=\varphi,\quad
\partial_t u(0)=\psi.
\end{equation}
We assume that 
the $k$-th component of the vector function $F$ takes the form 
$F^k(\partial u,\partial^2 u)=G^k(u,u,u)+H^k(u,u,u)$, where
\begin{eqnarray}
& &
G^k(u,v,w)
=
G_{ijl}^{k,\alpha\beta\gamma\delta}
\partial_\alpha u^i\partial_\beta v^j\partial_\gamma\partial_\delta w^l,\\
& &
H^k(u,v,w)=H_{ijl}^{k,\alpha\beta\gamma}
\partial_\alpha u^i\partial_\beta v^j\partial_\gamma w^l\nonumber
\end{eqnarray}
for real constants $G_{ijl}^{k,\alpha\beta\gamma\delta}$ 
and 
$H_{ijl}^{k,\alpha\beta\gamma}$. 
We refer to a term as {\it non-resonant} 
if $(i,j,l)\ne(k,k,k)$ in its coefficient. 
The remaining ones are said to be {\it resonant}.

Since our proof is based on the energy integral method, 
we naturally suppose the symmetry condition
\begin{equation}
G^{k,\alpha\beta\gamma\delta}_{ijl}
=
G^{k,\alpha\beta\delta\gamma}_{ijl}
=
G^{l,\alpha\beta\gamma\delta}_{ijk}.
\end{equation}
We are now in a position to recall the null condition 
in the setting of multiple speeds 
which Agemi and Yokoyama $\cite{AY}$ proposed: 
For every $k=1,\dots,m$ there holds
\begin{equation}
G_{kkk}^{k,\alpha\beta\gamma\delta}
X_\alpha X_\beta X_\gamma X_\delta
=
H_{kkk}^{k,\alpha\beta\gamma}
X_\alpha X_\beta X_\gamma
=0
\end{equation}
for all $X=(X_0,X_1,X_2)\in\{\,X\in{\mathbb R}^{1+2}\,:\,
X_0^2-c_k^2(X_1^2+X_2^2)=0\,\}$. 
The main theorem of this paper reads as follows.
\smallskip

\noindent{\bf Main Theorem}. {\it Assume 
the different-speed condition $(2.1)$, 
the symmetry condition $(3.4)$, 
and the null condition $(3.5)$. 
Let $\kappa\geq 9$. 
Then there exist positive constants 
$\varepsilon$ and $A$ with the following property$:$ 
If initial data is small so that
\begin{equation}
E_{\kappa-2}^{1/2}(u(0))
\exp
\bigl(
AE_\kappa^{1/2}(u(0))
(E_\kappa^{1/2}(u(0))+N_{\kappa-2}(u(0)))
\bigr)
<\varepsilon
\end{equation}
may hold, then the problem $(3.1)$-$(3.2)$ has a unique global in time 
solution satisfying 
\begin{eqnarray}
& &
E_\kappa(u(t))
\leq
4E_\kappa(u(0))(1+t)^{C\varepsilon^2},
\,\,\,
E_{\kappa-2}(u(t))
<4\varepsilon^2,\\
& &
\sum_{|a|\leq\kappa-3}
\|\langle x\rangle^{-1}\partial\Gamma^a u\|_{L^2((0,\,T)\times{\mathbb R}^2)}
\leq
C\varepsilon\log(2+T),\\
& &
N_{\kappa-2}(u(t))
\leq
N_{\kappa-2}(u(0))
+C\varepsilon^2
E_\kappa^{1/2}(u(0))
(1+t)^{(1/4)+C\varepsilon^2}
\end{eqnarray}
for all $t,\,T>0$.
}

{\it Remark.} The quantities 
$E_\kappa(u(0))$, $E_{\kappa-2}(u(0))$, and $N_{\kappa-2}(u(0))$ depend on 
the size of the initial data $(\varphi,\psi)$. 
Indeed, for given data $(\varphi,\psi)$, 
we can calculate the derivatives of the solution $u$ 
at $t=0$ up to the $\kappa$-th order 
by using the equation (3.1). 
In this way, we can determine these three quantities explicitly.
\section{Preliminaries}
\setcounter{equation}{0}
In this section, we collect several lemmas concerning 
commutation relations, some estimates of the null forms, and 
the Sobolev-type inequalities. 

We begin with the commutation relations. Let $[\cdot,\cdot]$ be 
the commutator. In addition to the well-known facts
\begin{equation}
[\partial_\alpha,\square]=0,\quad
[\Omega,\square]=0,\,\,\,\mbox{and}\quad
[S,\square]=-2\square,
\end{equation}
we need the commutation relations of the vector fields 
$\Gamma$ with respect to the nonlinear terms. 
Recall the nonlinear terms 
$G=(G^1,\dots,G^m)$ and $H=(H^1,\dots,H^m)$ defined in (3.3). 
Part (i) of the following lemma implies 
that the null structure is preserved upon differentiation, 
and Part (ii) together with (4.1) inductively shows that, 
for any $a$, the nonlinear term of the equation (4.4) 
also possesses the null structure. 
\smallskip

\noindent{\bf Lemma 4.1} (i) {\it 
For any $\Gamma^a$, the following equalities hold$:$
\begin{eqnarray}
& &
\Gamma^a G(u,v,w)
=
\sum_{b+c+d+e=a}
G_e(\Gamma^b u,\Gamma^c v,\Gamma^d w),\\
& &
\Gamma^a H(u,v,w)
=
\sum_{b+c+d+e=a}
H_e(\Gamma^b u,\Gamma^c v,\Gamma^d w).
\end{eqnarray}
Here each $G_e$ $($resp.\,$H_e)$ 
is a cubic nonlinear term of the form which 
$G$ $($resp.\,$H)$ has in $(3.3)$. In particular, 
$G_e=G$, $H_e=H$ if 
$b+c+d=a$ in $(4.2)$-$(4.3)$. 
Moreover, if the original 
nonlinearities $G$ and $H$ have the null structure $(3.5)$, 
then so does each of new nonlinearities $G_e$ and $H_e$. 

$({\rm ii})$ Let $u$ be a smooth solution of $(3.1)$-$(3.3)$. 
Then, for any $\Gamma^a$, the equalities
\begin{eqnarray}
& &
\square\Gamma^a u
=
\sum_{b+c+d+e=a}
G_e(\Gamma^b u,\Gamma^c u,\Gamma^d u)\\
& &
\hspace{1.1cm}
+
\sum_{b+c+d+e=a}
H_e(\Gamma^b u,\Gamma^c u,\Gamma^d u)
-
[\Gamma^a,\square]u\nonumber
\end{eqnarray}
hold.
}

{\it Proof}. See Lemma 4.1 of Sideris and Tu $\cite{ST}$. 
$\hfill\square$

The next lemma, which crucially comes into play 
in the estimates of lower-order energies, 
is the statement of gain of additional decay 
in nonlinearities with the null structure (3.5). 
\smallskip

\noindent{\bf Lemma 4.2} {\it For any smooth 
scalar functions $u$, $v$, $w$ and $z$, the following inequalities 
hold for $r\geq c_k t/2$$:$
\begin{eqnarray}
& &
|G_{kkk}^{k,\alpha\beta\gamma\delta}
\partial_\alpha u\partial_\beta v
\partial_\gamma\partial_\delta w|\\
& &
\leq
C\langle t\rangle^{-1}
\bigl[
|\Gamma u\|\partial v\|\partial^2 w|
+
|\partial u\|\Gamma v\|\partial^2 w|\nonumber\\
& &
\hspace{2.0cm}
+
|\partial u\|\partial v\|\partial\Gamma w|
+
\langle c_k t-r\rangle
|\partial u\|\partial v\|\partial^2 w|
\bigr],\nonumber\\
& &
|G_{kkk}^{k,\alpha\beta\gamma\delta}
\partial_\alpha u\partial_\beta v
\partial_\gamma w\partial_\delta z|\\
& &
\leq
C\langle t\rangle^{-1}
\bigl[
|\Gamma u\|\partial v\|\partial w\|\partial z|
+
|\partial u\|\Gamma v\|\partial w\|\partial z|
+
|\partial u\|\partial v\|\Gamma w\|\partial z|\nonumber\\
& &
\hspace{2.0cm}
+
|\partial u\|\partial v\|\partial w\|\Gamma z|
+
\langle c_k t-r\rangle
|\partial u\|\partial v\|\partial w\|\partial z|
\bigr],\nonumber\\
& &
|G_{kkk}^{k,\alpha\beta\gamma\delta}
\partial_\alpha \partial_\gamma u
\partial_\beta v
\partial_\delta w|\\
& &
\leq
C\langle t\rangle^{-1}
\bigl[
|\partial\Gamma u\|\partial v\|\partial w|
+
|\partial^2 u\|\Gamma v\|\partial w|\nonumber\\
& &
\hspace{2.0cm}
+
|\partial^2 u\|\partial v\|\Gamma w|
+
\langle c_k t-r\rangle|\partial^2 u\|\partial v\|\partial w|
\bigr],\nonumber\\
& &
|G_{kkk}^{k,\alpha\beta\gamma\delta}
\partial_\alpha u
\partial_\gamma\partial_\beta v
\partial_\delta w|\\
& &
\leq
C\langle t\rangle^{-1}
\bigl[
|\Gamma u\|\partial^2 v\|\partial w|
+
|\partial u\|\partial\Gamma v\|\partial w|\nonumber\\
& &
\hspace{2.0cm}
+
|\partial u\|\partial^2 v\|\Gamma w|
+
\langle c_k t-r\rangle|\partial u\|\partial^2 v\|\partial w|
\bigr],\nonumber\\
& &
|H_{kkk}^{k,\alpha\beta\gamma}
\partial_\alpha u\partial_\beta v\partial_\gamma w|\\
& &
\leq
C\langle t\rangle^{-1}
\bigl[
|\Gamma u\|\partial v\|\partial w|
+
|\partial u\|\Gamma v\|\partial w|\nonumber\\
& &
\hspace{2.0cm}
+
|\partial u\|\partial v\|\Gamma w|
+
\langle c_k t-r\rangle
|\partial u\|\partial v\|\partial w|
\bigr].\nonumber
\end{eqnarray}
}

\noindent{\it Proof}. We have only to mimic 
the proof of Lemma 5.1 of Sideris and Tu $\cite{ST}$. $\hfill\square$

The following lemma is concerned with Sobolev-type inequalities. 
\smallskip

\noindent{\bf Lemma 4.3} {\it The following inequalities hold for 
any smooth vector-valued 
function $u:(0,\infty)\times{\mathbb R}^2\to{\Bbb R}^m$, provided that 
the norms on the right-hand side are finite$:$
\begin{eqnarray}
& &
\langle r\rangle^{1/2}
|\partial u(t,x)|
\leq
CE_3^{1/2}(u(t)),\\
& &
\langle r\rangle^{1/2}\langle c_jt-r\rangle^{1/2}
|\partial u^j(t,x)|
\leq
CE_2^{1/4}(u(t))M_3^{1/2}(u(t)),\\
& &
\langle r\rangle^{1/2}
\langle c_j t-r\rangle
|\partial^2 u^j(t,x)|
\leq
CM_4(u(t)).
\end{eqnarray}
Moreover, for any $p$ with $2<p<\infty$, 
there exists a constant $C$ depending only on $p$ such that 
the inequality 
\begin{equation}
\langle t\rangle^{(1/2)-(1/p)}
\|\partial u(t)\|_{L^p({\mathbb R}^2)}
\leq
CE_2^{1/4}(u(t))M_3^{1/2}(u(t))
\end{equation}
holds.
}

{\it Proof}. The first three inequalities are proved in 
Lemma 1 of Sideris $\cite{Tom2}$. 
For the proof of (4.13), we have only to employ (4.11) and 
follow the proof of Lemma 2.2 (ii) of Katayama \cite{Ka}.
$\hfill\square$

{\it Remark}. The right-hand side of (4.11) takes 
the ``multiplicative'' form, 
which plays an important role in the proof of Lemma 5.3 below. 
\section{Weighted $L^2$-estimates}
\setcounter{equation}{0}
It is necessary to bound the weighted $L^2$-norm 
$M_\kappa(u(t))$ by $E_\kappa^{1/2}(u(t))$ 
for the completion of the energy integral argument. 
We carry out this by starting with the next crucial inequality 
due to Klainerman and Sideris $\cite{KS}$, estimating the nonlinear terms 
carefully, and doing a bootstrap argument. 
\smallskip

\noindent{\bf Lemma 5.1 (Klainerman--Sideris inequality)} 
{\it Let $\kappa\geq 2$. The inequality
\begin{equation}
M_\kappa(u(t))
\leq
C
\bigl(
E_\kappa^{1/2}(u(t))
+
\sum_{|a|\leq\kappa-2}
\|(t+r)\square\Gamma^a u(t)\|_{L^2({\mathbb R}^2)}
\bigr)
\end{equation}
holds for any smooth function $u$ with the finite norms 
on the right-hand side.
}

{\it Proof}. See Lemma 3.1 of Klainerman and Sideris $\cite{KS}$ 
and Lemma 7.1 of Sideris and Tu $\cite{ST}$. 
Note that their proof is obviously valid for 
$n=2$ as well as $n=3$.
$\hfill\square$

In the following, we denote by $[x]$ the greatest integer not greater than $x$.
\smallskip

\noindent{\bf Lemma 5.2} {\it 
Let $u$ be a smooth solution of $(3.1)$-$(3.2)$. 
Set $\kappa'=[(\kappa-1)/2]+3$. Then for all 
$|a|\leq \kappa-2$, it holds that
\begin{eqnarray}
& &
\|(t+r)\square\Gamma^a u(t)\|_{L^2({\mathbb R}^2)}\\
& &
\leq
C
\bigl(
E_{\kappa'}^{1/4}(u(t))M_{\kappa'}^{1/2}(u(t))
\bigr)^2
E_\kappa^{1/2}(u(t))\nonumber\\
& &
+
CE_{\kappa'}(u(t))M_{\kappa}(u(t))
+
CE_{\kappa'}^{1/2}(u(t))
E_{\kappa}^{1/2}(u(t))
M_{\kappa'}(u(t)).\nonumber
\end{eqnarray}
}

{\it Proof}. We may focus on the estimate of the 
$L^2$-norm of $t\square\Gamma^a u(t)$ because 
that of $r\square\Gamma^a u(t)$ is treated in a similar 
(in fact, easier) way. Set 
$p=[(\kappa-1)/2]$, so that $p+3=\kappa'$. 
It immediately follows from (4.4) that 
\begin{eqnarray}
& &
t
\|\square\Gamma^a u(t)\|_{L^2}
\leq
C\sum_{i,j,l}\sum_{{|b|+|c|+|d|}\atop{\leq\kappa-2}}
t
\bigl(
\|\partial\Gamma^b u^i(t)
\partial\Gamma^c u^j(t)
\partial^2\Gamma^d u^l(t)\|_{L^2}
\\
& &
\hspace{5.5cm}
+\|\partial\Gamma^b u^i(t)
\partial\Gamma^c u^j(t)
\partial\Gamma^d u^l(t)\|_{L^2}
\bigr).\nonumber
\end{eqnarray}
For the estimate of the second term on the right-hand side of (5.3), 
we may suppose $|b|+|c|\leq p$ without loss of generality. 
We get
\begin{eqnarray}
& &
\|\partial\Gamma^b u^i(t)
\partial\Gamma^c u^j(t)
\partial\Gamma^d u^l(t)\|_{L^2}\\
& &
\leq
\langle t\rangle^{-1}
\|\langle r\rangle^{1/2}
\langle
c_i t-r
\rangle^{1/2}
\partial\Gamma^b u^i(t)\|_{L^\infty}\nonumber\\
& &
\hspace{1cm}
\times
\|\langle r\rangle^{1/2}\langle c_j t-r\rangle^{1/2}
\partial\Gamma^c u^j(t)\|_{L^\infty}
\|\partial\Gamma^d u^l(t)\|_{L^2}\nonumber\\
& &
\leq
\langle t\rangle^{-1}
C
\bigl(
E_{\kappa'}^{1/4}(u(t))M_{\kappa'}^{1/2}(u(t))
\bigr)^2
E_\kappa^{1/2}(u(t))\nonumber
\end{eqnarray}
by using (4.11). 
For the first terms on the right-hand side of (5.3), 
we sort them out into two groups: $|b|+|c|\leq p$ or $|d|\leq p-1$. 
The first group is estimated as 
\begin{eqnarray}
& &
\dots
\leq
\langle t\rangle^{-1}
\|\langle r\rangle^{1/2}
\partial\Gamma^b u^i(t)\|_{L^\infty}
\|\langle r\rangle^{1/2}\partial\Gamma^c u^j(t)\|_{L^\infty}
\|\langle c_l t-r\rangle\partial^2 \Gamma^d u^l(t)\|_{L^2}
\\
& &
\hspace{0.6cm}
\leq
C\langle t\rangle^{-1}
E_{\kappa'}(u(t))M_\kappa(u(t))\nonumber
\end{eqnarray}
by (4.10). Otherwise, assuming $|b|\leq p$ as well as 
$|d|\leq p-1$ without loss of generality, we get
\begin{eqnarray}
& &
\dots\leq
\langle t\rangle^{-1}
\|\langle r\rangle^{1/2}\partial\Gamma^b u^i(t)\|_{L^\infty}
\|\partial\Gamma^c u^j(t)\|_{L^2}
\|\langle r\rangle^{1/2}\langle c_l t-r\rangle
\partial^2\Gamma^d u^l(t)\|_{L^\infty}\\
& &
\hspace{0.6cm}
\leq
C\langle t\rangle^{-1}
E_{\kappa'}^{1/2}(u(t))
M_{\kappa'}(u(t))
E_{\kappa-1}^{1/2}(u(t))\nonumber
\end{eqnarray}
by (4.10), (4.12), which completes the proof of (5.2).
$\hfill\square$
\smallskip

\noindent{\bf Lemma 5.3} {\it 
Let $\kappa\geq 9$, $\mu=\kappa-2$. 
There exists a small, positive constant 
$\varepsilon_0$ with the following property$:$ 
Suppose that, 
for a local smooth solution $u$ of $(3.1)$-$(3.2)$, 
the supremum of $E_\mu^{1/2}(u(t))$ over an interval 
$[0,T)$ is sufficiently small so that 
\begin{equation}
\sup_{0\leq t<T}E_\mu^{1/2}(u(t))
\leq
\varepsilon_0
\end{equation}
may hold. 
Then
\begin{equation}
M_\mu(u(t))
\leq
CE_\mu^{1/2}(u(t)),\,\,0\leq t<T
\end{equation}
and 
\begin{equation}
M_\kappa(u(t))
\leq
CE_\kappa^{1/2}(u(t)),\,\,0\leq t<T
\end{equation}
hold with a constant $C$ independent of $T$.
}

{\it Remark}. This lemma is actually valid for $\kappa\geq 8$. 
We have assumed $\kappa\geq 9$ for the latter use. 

{\it Proof}. Set $\mu'=[(\mu-1)/2]+3$. 
Denoting by $\delta$ the supremum of $E_\mu^{1/2}(u(t))$ over 
the interval $[0,T)$, 
we see that Lemma 5.1 and Lemma 5.2 imply for $0\leq t<T$
\begin{eqnarray}
& &
M_\mu(u(t))
\leq
CE_\mu^{1/2}(u(t))
+
C
\bigl(
E_{\mu'}^{1/4}(u(t))
M_{\mu'}^{1/2}(u(t))
\bigr)^2
E_\mu^{1/2}(u(t))\\
& &
\hspace{2cm}
+
CE_{\mu'}(u(t))M_\mu(u(t))
+
CE_{\mu'}^{1/2}(u(t))E_\mu^{1/2}(u(t))M_{\mu'}(u(t))\nonumber\\
& &
\hspace{1.7cm}
\leq
CE_\mu^{1/2}(u(t))
+
C
\bigl(
\delta^{1/2}M_{\mu'}^{1/2}(u(t))
\bigr)^2
E_\mu^{1/2}(u(t))\nonumber\\
& &
\hspace{2cm}
+
C\delta^2M_\mu(u(t))
+
C\delta^2M_{\mu'}(u(t))\nonumber\\
& &
\hspace{1.7cm}
\leq
CE_\mu^{1/2}(u(t))
+
C\delta^2M_\mu(u(t)),\nonumber
\end{eqnarray}
which immediately yields (5.8) 
if $\delta$ is sufficiently small. 

As for (5.9), we first note that the inequality 
$\kappa':=[(\kappa-1)/2]+3\leq\mu$ holds. 
Proceeding as in (5.10) and using (5.8), 
we easily see that
\begin{equation}
M_\kappa(u(t))
\leq
CE_\kappa^{1/2}(u(t))
+
C\delta^2M_\kappa(u(t))
+
CE_\kappa^{1/2}(u(t)),
\end{equation}
which yields (5.9).
$\hfill\square$
\section{Energy estimates}

\setcounter{equation}{0}
Following the strategy in Sideris $\cite{Tom3}$ 
and Sideris and Tu $\cite{ST}$, 
we accomplish the energy integral argument 
by deriving a pair of coupled differential inequalities 
for a higher-order 
energy $E_\kappa (u(t))$, $\kappa\geq 9$ and 
a lower-order energy $E_\mu (u(t))$, $\mu=\kappa-2$. 
Since the equation is quasi-linear, 
we must actually consider modified energies which 
are equivalent to the original 
ones for small solutions. 

For initial data $(\varphi,\psi)$, 
let us assume 
$E_\mu^{1/2}(u(0))<
\varepsilon$ for a sufficiently 
small $\varepsilon>0$ such that $2\varepsilon\leq\varepsilon_0$ 
(see (5.7) for $\varepsilon_0$). 
By the standard local existence theorem, 
we know that a unique smooth solution exists 
locally in time. 
Suppose that $T_0$ is the supremum of all $T>0$ for which 
$E_\mu^{1/2}(u(t))<2\varepsilon$ for all $0\leq t<T$. 
It is shown that $E_\mu^{1/2}(u(t))<2\varepsilon$ on the closed interval 
$0\leq t\leq T_0$, therefore we can continue the local solution 
to all time.

Suppose $0\leq t<T_0$ in what follows. Denoting by 
$\langle\cdot,\cdot\rangle$ the scalar product in ${\Bbb R}^m$, 
we have for each $n=1,\dots,\kappa$ $(\kappa\geq 9)$
\begin{eqnarray}
& &
E'_n(u(t))
=
\sum_{|a|\leq n-1}
\int_{{\Bbb R}^2}
\langle
\square\Gamma^a u(t),\partial_t\Gamma^a u(t)
\rangle dx\\
& &
\hspace{1.66cm}
=
\sum_{{1\leq k\leq m}\atop{|a|=n-1}}
\int_{{\Bbb R}^2}
G_{ijl}^{k,\alpha\beta\gamma\delta}
\partial_\alpha u^i
\partial_\beta u^j
\partial_\gamma\partial_\delta\Gamma^a u^l
\partial_t\Gamma^a u^k dx\nonumber\\
& &
\hspace{2cm}
+
\sum_{{b+c+d+e=a}\atop{|a|\leq n-1,d\ne a}}
\int_{{\Bbb R}^2}
\langle
G_e(\Gamma^b u,\Gamma^c u,\Gamma^d u),\partial_t\Gamma^a u
\rangle
dx\nonumber\\
& &
\hspace{2cm}
+
\sum_{b+c+d+e=a}
\int_{{\Bbb R}^2}
\langle
H_e(\Gamma^b u,\Gamma^c u,\Gamma^d u),\partial_t\Gamma^a u
\rangle
dx\nonumber\\
& &
\hspace{2cm}
-\int_{{\Bbb R}^2}
\langle
[\Gamma^a,\square]u,\partial_t\Gamma^a u
\rangle
dx.\nonumber
\end{eqnarray}
The loss of derivatives which has occurred in the first term 
on the right-hand side is prevented by the symmetry condition (3.4) 
as follows:
\begin{eqnarray}
& &
\sum_{k=1}^m
\int_{{\Bbb R}^2}
G_{ijl}^{k,\alpha\beta\gamma\delta}
\partial_\alpha u^i
\partial_\beta u^j
\partial_\gamma\partial_\delta\Gamma^a u^l
\partial_t\Gamma^a u^k dx\\
& &
=
\sum_{k=1}^m
\int_{{\Bbb R}^2}
G_{ijl}^{k,\alpha\beta\gamma\delta}
\partial_\gamma
\bigl(
\partial_\alpha u^i
\partial_\beta u^j
\partial_\delta\Gamma^a u^l
\partial_t\Gamma^a u^k
\bigr)dx\nonumber\\
& &
\hspace{0.5cm}
-
\int_{{\Bbb R}^2}
G_{ijl}^{k,\alpha\beta\gamma\delta}
\bigl[
\partial_\gamma
\bigl(
\partial_\alpha u^i
\partial_\beta u^j
\bigr)
\partial_\delta\Gamma^a u^l
\partial_t\Gamma^a u^k\nonumber\\
& &
\hspace{3cm}
+
\partial_\alpha u^i
\partial_\beta u^j
\partial_\delta \Gamma^a u^l
\partial_t\partial_\gamma \Gamma^a u^k
\Bigr]dx\nonumber\\
& &
=
\partial_t
\sum_{k=1}^m
\int_{{\Bbb R}^2}
G_{ijl}^{k,\alpha\beta 0\delta}
\partial_\alpha u^i
\partial_\beta u^j
\partial_\delta\Gamma^a u^l
\partial_t\Gamma^a u^k dx\nonumber\\
& &
\hspace{0.5cm}
-
\sum_{k=1}^m
\int_{{\Bbb R}^2}
G_{ijl}^{k,\alpha\beta\gamma\delta}
\partial_\gamma
\bigl(
\partial_\alpha u^i
\partial_\beta u^j
\bigr)
\partial_\delta\Gamma^a u^l
\partial_t\Gamma^a u^k dx\nonumber\\
& &
\hspace{0.5cm}
-
\sum_{k=1}^m
\int_{{\Bbb R}^2}
\frac12
G_{ijl}^{k,\alpha\beta\gamma\delta}
\partial_\alpha u^i
\partial_\beta u^j
\partial_t
\bigl(
\partial_\delta\Gamma^a u^l
\partial_\gamma\Gamma^a u^k
\bigr)dx\nonumber\\
& &
=
\partial_t
\sum_{k=1}^m
\int_{{\Bbb R}^2}
\frac12
G_{ijl}^{k,\alpha\beta\lambda\delta}
\eta^{\gamma}_{\lambda}
\partial_\alpha u^i
\partial_\beta u^j
\partial_\delta\Gamma^a u^l
\partial_\gamma\Gamma^a u^k dx\nonumber\\
& &
\hspace{0.5cm}
-
\sum_{k=1}^m
\int_{{\Bbb R}^2}
G_{ijl}^{k,\alpha\beta\gamma\delta}
\partial_\gamma
\bigl(
\partial_\alpha u^i
\partial_\beta u^j
\bigr)
\partial_\delta\Gamma^a u^l
\partial_t\Gamma^a u^k dx\nonumber\\
& &
\hspace{0.5cm}
+
\sum_{k=1}^m
\int_{{\Bbb R}^2}
\frac12
G_{ijl}^{k,\alpha\beta\gamma\delta}
\partial_t
\bigl(
\partial_\alpha u^i
\partial_\beta u^j
\bigr)
\partial_\delta \Gamma^a u^l
\partial_\gamma\Gamma^a u^k dx.\nonumber
\end{eqnarray}
Here 
$\eta^\gamma_\lambda:=\mbox{diag}(1,-1,-1)$. 
Therefore, introducing the modified energy
\begin{eqnarray}
& &
{\tilde E}_n(u(t)):=E_n(u(t))\\
& &
\hspace{2cm}
-
\sum_{{|a|=n-1}\atop{1\leq k\leq m}}
\int_{{\Bbb R}^2}
\frac12
G_{ijl}^{k,\alpha\beta\lambda\delta}
\eta_\lambda^\gamma
\partial_\alpha u^i
\partial_\beta u^j
\partial_\delta\Gamma^a u^l
\partial_\gamma\Gamma^a u^k dx,\nonumber
\end{eqnarray}
we finally have
\begin{eqnarray}
& &
{\tilde E}'_n(u(t))
=
\sum_{{b+c+d+e=a}\atop{|a|\leq n-1,d\ne a}}
\int_{{\Bbb R}^2}
\langle
G_e(\Gamma^b u,
\Gamma^c u,\Gamma^d u),\partial_t\Gamma^a u
\rangle dx\\
& &
\hspace{1.7cm}
-
\sum_{{|a|=n-1}\atop{1\leq k\leq m}}
\int_{{\Bbb R}^2}
G_{ijl}^{k,\alpha\beta\gamma\delta}
\partial_\gamma
\bigl(
\partial_\alpha u^i
\partial_\beta u^j
\bigr)
\partial_\delta \Gamma^a u^l
\partial_t\Gamma^a u^k dx\nonumber\\
& &
\hspace{1.7cm}
+
\sum_{{|a|=n-1}\atop{1\leq k\leq m}}
\frac12
\int_{{\Bbb R}^2}
G_{ijl}^{k,\alpha\beta\gamma\delta}
\partial_t
\bigl(
\partial_\alpha u^i
\partial_\beta u^j
\bigr)
\partial_\delta\Gamma^a u^l
\partial_\gamma\Gamma^a u^k dx\nonumber\\
& &
\hspace{1.7cm}
+\sum_{b+c+d+e=a}
\int_{{\Bbb R}^2}
\langle
H_e(\Gamma^b u,\Gamma^c u,\Gamma^d u),
\partial_t\Gamma^a u
\rangle dx\nonumber\\
& &
\hspace{1.7cm}
-\int_{{\Bbb R}^2}
\langle
[\Gamma^a,\square]u,\partial_t\Gamma^a u
\rangle dx.\nonumber
\end{eqnarray}
We also note that, under the smallness 
of $E_\mu^{1/2}(u(t))$ 
$(0\leq t<T_0)$ with $\mu=\kappa-2$, the inequality
\begin{equation}
\frac12
E_n(u(t))
\leq
{\tilde E}_n(u(t))
\leq
2E_n(u(t)),
\,\,\,
n=1,\dots,\kappa
\end{equation}
holds by the Sobolev embedding.

We plan our energy integral method, 
allowing the higher-order energy 
$E_\kappa(u(t))$ $(\kappa\geq 9)$ to grow polynomially in time 
but bounding the lower-order energy 
$E_\mu(u(t))$ $(\mu=\kappa-2)$ uniformly in time. 
(See (3.7) above.) 
Let us start with the estimate of the higher-order energy. 
Setting $n=\kappa$ in (6.4), we have
\begin{eqnarray}
& &
{\tilde E}'_{\kappa}(u(t))
\leq
\sum_{i,j,l}
\sum_{|a|\leq \kappa-1}
\sum_{{|b|+|c|+|d|\leq |a|}\atop{d\ne a}}
\|\partial\Gamma^b u^i
\partial\Gamma^c u^j
\partial^2\Gamma^d u^l\|_{L^2}
\|\partial\Gamma^a u\|_{L^2}\\
& &
\hspace{1.7cm}
+
\sum_{i,j,l}
\sum_{|a|\leq\kappa-1}
\sum_{|b|+|c|+|d|\leq |a|}
\|\partial\Gamma^b u^i
\partial\Gamma^c u^j
\partial\Gamma^d u^l\|_{L^2}
\|\partial\Gamma^a u\|_{L^2}.\nonumber
\end{eqnarray}
Set $q=[\kappa/2]$. 
Note that $q+3\leq\mu$ because of $\kappa\geq 9$. 
Supposing $|b|+|c|\leq q$ without loss of generality, 
we bound the second term as 
\begin{eqnarray}
& &
\|\partial\Gamma^b u^i
\partial\Gamma^c u^j
\partial\Gamma^d u^l\|_{L^2}\\
& &
\leq
C
\langle t\rangle^{-1}
\|\langle r\rangle^{1/2}
\langle c_i t-r\rangle^{1/2}
\partial\Gamma^b u^i\|_{L^\infty}
\|\langle r\rangle^{1/2}
\langle c_j t-r\rangle^{1/2}\partial\Gamma^c u^j\|_{L^\infty}
\|\partial\Gamma^d u^l\|_{L^2}\nonumber\\
& &
\leq
C\langle t\rangle^{-1}
\bigl(
E_\mu^{1/4}(u(t))M_\mu^{1/2}(u(t))
\bigr)^2
E_\kappa^{1/2}(u(t))
\leq
C\langle t\rangle^{-1}
E_\mu(u(t))E_{\kappa}^{1/2}(u(t)).\nonumber
\end{eqnarray}
Here we have employed (5.8) at the third inequality. 
As for the first terms on the right-hand side of (6.6), 
we sort them out into two groups: $|b|+|c|\leq q$ or 
$|d|\leq q-1$. The first group is estimated as in 
(5.5) and (6.7):
\begin{eqnarray}
& &
\|\partial\Gamma^b u^i
\partial\Gamma^c u^j
\partial^2\Gamma^d u^l\|_{L^2}
\leq
C\langle t\rangle^{-1}
\bigl(
E_\mu^{1/4}(u(t))M_\mu^{1/2}(u(t))
\bigr)^2M_\kappa(u(t))\\
& &
\hspace{4.14cm}
\leq
C\langle t\rangle^{-1}
E_\mu(u(t))E_\kappa^{1/2}(u(t)).\nonumber
\end{eqnarray}
Otherwise, assuming $|b|\leq q$ in addition to 
$|d|\leq q-1$ without loss of generality, we get as in (5.6)
\begin{eqnarray}
& &
\|\partial\Gamma^b u^i
\partial\Gamma^c u^j
\partial^2\Gamma^d u^l\|_{L^2}
\leq
C\langle t\rangle^{-1}
E_\mu^{1/2}(u(t))M_\mu(u(t))E_{\kappa}^{1/2}(u(t))\\
& &
\hspace{4.14cm}
\leq
C\langle t\rangle^{-1}E_\mu(u(t))E_\kappa^{1/2}(u(t)).\nonumber
\end{eqnarray}
Taking account of the equivalence between 
$E_n$ and ${\tilde E}_n$, we get from 
(6.6)-(6.9)
\begin{equation}
{\tilde E}'_\kappa(u(t))
\leq
C\langle t\rangle^{-1}
{\tilde E}_\mu(u(t))
{\tilde E}_\kappa(u(t)).
\end{equation}
{\bf Lower-order Energy}. The crucial part in the proof of 
global existence is to bound the lower-order energy 
$E_\mu(u(t))$ $(\mu=\kappa-2)$ uniformly in time. 
For the purpose, 
we exploit the difference of propagation speeds 
as well as an improved decay rate of solutions inside the cone 
to sharpen the decay estimates presented above, 
when $|a|\leq\mu$. 
Moreover, the space-time $L^2((0,\infty)\times{\mathbb R}^2)$ estimate and 
the ${\dot H}^{1/2}({\mathbb R}^2)$ estimate 
play an auxiliary role. 

Set $c_0:=\min\{c_i/2:i=1,\dots,m\}$ and 
$\mu=\kappa-2$ $(\kappa\geq 9)$. 
Setting $n=\mu$ in (6.4), we estimate the resulting terms 
on the right-hand side. 
Divide the integral region ${\Bbb R}^2$ into two parts: 
inside the cone 
$\{(t,x):|x|\leq c_0 t\}$ and 
away from the spatial origin 
$\{(t,x):|x|\geq c_0 t\}$. 
\smallskip

\noindent{\it Inside the cone}. Here we exploit an improved decay rate 
of solutions. 
The space-time $L^2((0,\infty)\times{\mathbb R}^2)$ estimate 
also comes into play. 
The contribution from the quasi-linear terms is bounded by
\begin{equation}
\sum_{i,j,l}
\sum_{|a|\leq \mu-1}
\sum_{{|b|+|c|+|d|\leq |a|}\atop{d\ne a}}
\|\partial\Gamma^b u^i
\partial\Gamma^c u^j
\partial^2\Gamma^d u^l\|_{L^2(r<c_0t)}
\|\partial\Gamma^a u\|_{L^2}.
\end{equation}
We may suppose $|b|\leq [\mu/2]$ 
without loss of generality. 
It then follows from (4.11) and (5.8) that 
\begin{eqnarray}
& &
\|\partial\Gamma^b u^i
\partial\Gamma^c u^j
\partial^2\Gamma^d u^l\|_{L^2}\\
& &
\leq
\langle t\rangle^{-3/2}
\|\langle c_i t-r\rangle^{1/2}
\partial\Gamma^b u^i\|_{L^\infty(r<c_0 t)}
\|\partial\Gamma^c u^j\|_{L^\infty}
\|\langle c_l t-r\rangle\partial^2\Gamma^d u^l\|_{L^2(r<c_0 t)}\nonumber\\
& &
\leq
C\langle t\rangle^{-3/2}
\bigl(
E_{|b|+2}^{1/4}(u(t))
M_{|b|+3}^{1/2}(u(t))
\bigr)
E_{|c|+3}^{1/2}(u(t))
M_\mu(u(t))\nonumber\\
& &
\leq
C\langle t\rangle^{-3/2}
E_\kappa^{1/2}(u(t))
E_\mu(u(t)),\nonumber
\end{eqnarray}
where we have used 
$|b|+3\leq[\mu/2]+3\leq\mu$, 
$|c|+3\leq\kappa$. 
Concerning the contribution from the semi-linear parts, 
we see from (6.4) that it is bounded by
\begin{eqnarray}
\sum_{i,j,l}
\sum_{|a|\leq \mu-1}
\sum_{|b|+|c|+|d|\leq |a|}
\|\partial\Gamma^b u^i
\partial\Gamma^c u^j
\partial\Gamma^d u^l\|_{L^2(r<c_0 t)}
\|\partial\Gamma^a u\|_{L^2}.
\end{eqnarray}
Assume $|b|+|c|\leq[\mu/2]$ without loss of generality. 
Proceeding quite differently from how we did in (6.23) of \cite{Hi4}, 
we get
\begin{eqnarray}
& &
\|\partial\Gamma^b u^i(t)
\partial\Gamma^c u^j(t)
\partial\Gamma^d u^l(t)\|_{L^2(r<c_0 t)}\\
& &
\leq
C\langle t\rangle^{-1}
\|\langle c_it-r\rangle^{1/2}\langle r\rangle^{1/2}
\partial\Gamma^b u^i(t)\|_{L^\infty(r<c_0 t)}\nonumber\\
& &
\hspace{1.3cm}
\times
\|\langle c_j t-r\rangle^{1/2}\langle r\rangle^{1/2}
\partial\Gamma^c u^j(t)\|_{L^\infty(r<c_0 t)}
\|\langle r\rangle^{-1}\partial\Gamma^d u^l(t)\|_{L^2(r<c_0 t)}\nonumber\\
& &
\leq
C\langle t\rangle^{-1}
\bigl(
E_{|b|+2}^{1/4}(u(t))M_{|b|+3}^{1/2}(u(t))
\bigr)
\bigl(
E_{|c|+2}^{1/4}(u(t))M_{|c|+3}^{1/2}(u(t))
\bigr)
\nonumber\\
& &
\hspace{1.3cm}
\times
\|\langle r\rangle^{-1}\partial\Gamma^d u^l(t)\|_{L^2({\mathbb R}^2)}
\nonumber\\
& &
\leq
C\langle t\rangle^{-1}
E_\mu(u(t))
\sum_{|d|\leq\mu-1}
\|\langle r\rangle^{-1}\partial\Gamma^d u(t)\|_{L^2({\mathbb R}^2)}.\nonumber
\end{eqnarray}
The estimate inside the cone has been finished.
\smallskip

\noindent{\it Away from the spatial origin}. Here the difference of 
propagation speeds comes into play. 
Moreover, we employ the null condition (3.5) 
for the estimates of resonance terms. 
\smallskip

\noindent{\it Non-resonance}. Let us start with non-resonance terms. 
Our task is to estimate the contribution from quasi-linear 
terms
\begin{equation}
\sum_{(i,j,l)\ne(k,k,k)}
\sum_{|a|\leq\mu-1}
\sum_{{|b|+|c|+|d|\leq|a|}\atop{d\ne a}}
\|\partial\Gamma^b u^i
\partial\Gamma^c u^j
\partial^2\Gamma^d u^l
\partial\Gamma^a u^k
\|_{L^1(r>c_0 t)}
\end{equation}
and the contribution from semi-linear terms
\begin{equation}
\sum_{(i,j,l)\ne(k,k,k)}
\sum_{|a|\leq\mu-1}
\sum_{|b|+|c|+|d|\leq |a|}
\|\partial\Gamma^b u^i
\partial\Gamma^c u^j
\partial\Gamma^d u^l
\partial_t\Gamma^a u^k\|_{L^1(r>c_0 t)}.
\end{equation}
In estimating the $L^1$-norm in (6.15) we separate 
two cases: $i=j=l$ or otherwise. 
In the former case, noting $i\ne k$, we have
\begin{eqnarray}
& &
\|\partial\Gamma^b u^i\partial\Gamma^c u^i
\partial^2\Gamma^d u^i\partial\Gamma^a u^k\|_{L^1(r>c_0t)}\\
& &
\leq
\langle t\rangle^{-3/2}
\|\langle r\rangle^{1/2}\partial\Gamma^b u^i\|_{L^\infty}
\|\partial\Gamma^c u^i\|_{L^2}\nonumber\\
& &
\hspace{3cm}
\times
\|\langle c_it-r\rangle\partial^2\Gamma^d u^i\|_{L^2}
\|\langle r\rangle^{1/2}\langle c_kt-r\rangle^{1/2}\partial\Gamma^a u^k\|
_{L^\infty}\nonumber\\
& &
\leq
C\langle t\rangle^{-3/2}
E_{|b|+3}^{1/2}(u(t))
E_{|c|+1}^{1/2}(u(t))
M_{|d|+2}(u(t))
\bigl(
E_{|a|+2}^{1/4}M_{|a|+3}^{1/2}(u(t))
\bigr)\nonumber\\
& &
\leq
C\langle t\rangle^{-3/2}
E_{\mu+2}^{1/2}(u(t))
E_\mu^{1/2}(u(t))
M_\mu(u(t))E_{\mu+2}^{1/2}(u(t))\nonumber\\
& &
\leq
C\langle t\rangle^{-3/2}E_\kappa(u(t))E_\mu(u(t)).\nonumber
\end{eqnarray}
Otherwise, it is easy to get 
\begin{eqnarray}
& &
\|\partial\Gamma^b u^i\partial\Gamma^c u^j
\partial^2\Gamma^d u^l\partial\Gamma^a u^k\|_{L^1(r>c_0 t)}\\
& &
\leq
\langle t\rangle^{-3/2}
\|\langle r\rangle^{1/2}\langle c_i t-r\rangle^{1/2}
\partial\Gamma^b u^i\|_{L^\infty}\nonumber\\
& &
\hspace{1cm}
\times
\|\langle r\rangle^{1/2}
\langle c_j t-r\rangle^{1/2}
\partial\Gamma^c u^j\|_{L^\infty}
\|\langle c_l t-r\rangle\partial^2\Gamma^d u^l\|_{L^2}
\|\partial\Gamma^a u^k\|_{L^2}\nonumber\\
& &
\leq
C\langle t\rangle^{-3/2}
E_\kappa(u(t))E_\mu(u(t)).\nonumber
\end{eqnarray}
As for (6.16), we may suppose $i\ne k$ without loss of generality. 
We obtain
\begin{eqnarray}
& &
\|\partial\Gamma^b u^i
\partial\Gamma^c u^j
\partial\Gamma^d u^l
\partial_t \Gamma^a u^k\|_{L^1(r>c_0 t)}\\
& &
\leq
\langle t\rangle^{-3/2}
\|\langle r\rangle^{1/2}
\langle c_i t-r\rangle^{1/2}
\partial\Gamma^b u^i\|_{L^\infty(r>c_0 t)}\nonumber\\
& &
\hspace{1cm}
\times
\|\partial\Gamma^c u^j\partial\Gamma^d u^l\|_{L^1}
\|\langle r\rangle^{1/2}
\langle c_k t-r\rangle^{1/2}
\partial_t \Gamma^a u^k\|_{L^\infty(r>c_0 t)}\nonumber\\
& &
\leq
C\langle t\rangle^{-3/2}
\bigl(
E_{|b|+2}^{1/4}(u(t))
M_{|b|+3}^{1/2}(u(t))
\bigr)
E_\mu(u(t))
\bigl(
E_{|a|+2}^{1/4}(u(t))
M_{|a|+3}^{1/2}(u(t))
\bigr)\nonumber\\
& &
\leq
C\langle t\rangle^{-3/2}
E_\kappa(u(t))
E_\mu(u(t)).\nonumber
\end{eqnarray}
Therefore the estimates of non-resonance terms away from 
the spatial origin have been completed. 
\smallskip

\noindent{\it Resonance}. The resonance terms remain to be estimated 
away from the spatial origin. It is just the place 
where the null condition comes into play. 
Without the null condition, the solution 
may become singular in finite time (see, e.g., Zhou and Han $\cite{ZH}$). 
In view of Lemma 4.2 and (6.4), the estimate is reduced to bounding
\begin{eqnarray}
& &
\sum_{1\leq k\leq m}
\sum_{|a|\leq \mu-1}
\sum_{{|b|+|c|+|d|\leq |a|}\atop{d\ne a}}
\langle t\rangle^{-1}
\bigl(
\|\Gamma^{b+1}u^k
\partial\Gamma^c u^k
\partial^2\Gamma^d u^k\|_{L^2(r>c_0 t)}\\
& &
\hspace{1cm}
+
\|\partial\Gamma^b u^k\partial\Gamma^c u^k
\partial\Gamma^{d+1} u^k\|_{L^2(r>c_0 t)}\nonumber\\
& &
\hspace{1cm}
+
\|\langle c_k t-r\rangle
\partial\Gamma^b u^k
\partial\Gamma^c u^k
\partial^2\Gamma^d u^k\|_{L^2(r>c_0 t)}
\bigr)
\|\partial\Gamma^a u\|_{L^2}\nonumber\\
& &
+
\sum_{1\leq k\leq m}
\sum_{|a|\leq \mu-1}
\sum_{|b|+|c|+|d|\leq |a|}
\langle t\rangle^{-1}
\bigl(
\|\Gamma^{b+1}u^k
\partial\Gamma^c u^k
\partial\Gamma^d u^k\|_{L^2(r>c_0 t)}\nonumber\\
& &
\hspace{1cm}
+
\|\langle c_k t-r\rangle\partial\Gamma^b u^k
\partial\Gamma^c u^k
\partial\Gamma^d u^k\|_{L^2(r>c_0 t)}
\bigr)
\|\partial\Gamma^a u\|_{L^2}.\nonumber
\end{eqnarray}
Here, by $b+1$, we mean any sequence of length $|b|+1$. 

We proceed differently from how we did in (6.40) of \cite{Hi4}. 
Using (4.10) and the Hardy inequality of order $1/2$, 
we estimate the first norm on the right-hand side 
of (6.20) as 
\begin{eqnarray}
& &
C\langle t\rangle^{-1/2}
\|
r^{-1/2}\Gamma^{b+1}u^k
r^{1/2}\partial\Gamma^c u^k
\langle r\rangle^{1/2}\partial^2\Gamma^d u^k
\|_{L^2(r>c_0 t)}\\
& &
\leq
C\langle t\rangle^{-1/2}
\|r^{-1/2}\Gamma^{b+1}u^k\|_{L^2}
\|r^{1/2}\partial\Gamma^c u^k\|_{L^\infty}
\|\langle r\rangle^{1/2}\partial^2\Gamma^d u^k\|_{L^\infty}\nonumber\\
& &
\leq
C\langle t\rangle^{-1/2}
\|\,|D_x|^{1/2}\Gamma^{b+1}u^k\,\|_{L^2}
E_{|c|+3}^{1/2}(u(t))
E_{|d|+3}^{1/2}(u(t))\nonumber\\
& &
\leq
C\langle t\rangle^{-1/2}
\|\,|D_x|^{1/2}\Gamma^{b+1}u^k\,\|_{L^2}
E_\mu^{1/2}(u(t))E_\kappa^{1/2}(u(t))
\nonumber\\
& &
\leq
C\langle t\rangle^{-1/2}
\bigl(
\sum_{|a|\leq\mu}
\|\Gamma^a u(t)\|_{{\dot H}^{1/2}}
\bigr)
E_\mu^{1/2}(u(t))E_\kappa^{1/2}(u(t)).\nonumber
\end{eqnarray}
Assuming $|b|\leq |c|$ without loss of generality, we estimate 
the second and third norms on the right-hand side of (6.20) as 
\begin{eqnarray}
& &
\langle t\rangle^{-1/2}
\|\langle r\rangle^{1/2}
\partial\Gamma^b u^k\|_{L^\infty(r>c_0 t)}
\|\partial\Gamma^c u^k\|_{L^\infty}
\|\partial\Gamma^{d+1}u^k\|_{L^2}\\
& &
+
\langle t\rangle^{-1/2}
\|\langle r\rangle^{1/2}
\partial\Gamma^b u^k\|_{L^\infty(r>c_0 t)}
\|\partial\Gamma^c u^k\|_{L^\infty}
\|\langle c_k t-r\rangle\partial^2\Gamma^d u^k\|_{L^2}\nonumber\\
& &
\leq
C\langle t\rangle^{-1/2}
E_\kappa^{1/2}(u(t))
E_\mu(u(t)).\nonumber
\end{eqnarray}
The remaining terms in (6.20) are estimated as 
\begin{eqnarray}
& &
\langle t\rangle^{-1/2}
\|r^{-1/2}\Gamma^{b+1}u^k
r^{1/2}\partial\Gamma^c u^k
\langle r\rangle^{1/2}
\partial\Gamma^d u^k\|_{L^2(r>c_0 t)}\\
& &
+\langle t\rangle^{-1}
\|\langle r\rangle^{1/2}\langle c_k t-r\rangle^{1/2}\partial\Gamma^b u^k
\langle r\rangle^{1/2}\langle c_k t-r\rangle^{1/2}\partial\Gamma^c u^k
\partial\Gamma^d u^k\|_{L^2(r>c_0 t)}\nonumber\\
& &
\leq
C\langle t\rangle^{-1/2}
\bigl(
\sum_{|a|\leq\mu}
\|\Gamma^a u(t)\|_{{\dot H}^{1/2}}
\bigr)
E_\mu^{1/2}(u(t))E_\kappa^{1/2}(u(t))\nonumber\\
& &
+
C\langle t\rangle^{-1}
E_\kappa^{1/2}(u(t))
E_\mu(u(t)),\nonumber
\end{eqnarray}
thanks to (4.10), (4.11), and the Hardy inequality of order $1/2$. 
The estimate of (6.20) has been finished. 

Collecting the estimates of ${\tilde E}'_\mu(u(t))$ and 
taking (6.5) into account, we have finally obtained
\begin{eqnarray}
& &
{\tilde E}'_\mu(u(t))
\leq
C\langle t\rangle^{-3/2}
{\tilde E}^{1/2}_\kappa(u(t)){\tilde E}_\mu^{3/2}(u(t))\\
& &
\hspace{1.7cm}
+
C\langle t\rangle^{-1}
{\tilde E}^{3/2}_\mu(u(t))
\sum_{|a|\leq\mu-1}
\|\langle r\rangle^{-1}\partial\Gamma^a u(t)\|_{L^2({\mathbb R}^2)}
\nonumber\\
& &
\hspace{1.7cm}
+C\langle t\rangle^{-3/2}
{\tilde E}_\kappa(u(t)){\tilde E}_\mu(u(t))\nonumber\\
& &
\hspace{1.7cm}
+
C\langle t\rangle^{-3/2}
{\tilde E}_\kappa^{1/2}(u(t))
{\tilde E}_\mu(u(t))
\bigl(
{\tilde E}_\mu^{1/2}(u(t))
+
\sum_{|a|\leq\mu}\|\Gamma^a u(t)\|_{{\dot H}^{1/2}}
\bigr).\nonumber
\end{eqnarray}

\noindent{\it Space-time $L^2$ estimate}. 
Two more ingredients are needed to accomplish 
the energy integration argument. 
One is the space-time $L^2$ estimate, 
and the other is the ${\dot H}^{1/2}({\mathbb R}^2)$ estimate. 
See the right-hand side of (6.24). 
For the former, we utilize the localized energy estimate 
of Smith and Sogge \cite{SS}: for $n\geq 1$ and 
$0\leq\gamma\leq (n-1)/2$, there holds
\begin{equation}
\|\beta(\exp(it|D_x|)g)\|_{L^2({\mathbb R};H^\gamma({\mathbb R}^n))}
\leq
C
\|\,|D_x|^\gamma g\,\|_{L^2({\mathbb R}^n)},
\end{equation}
here $\beta\in C_0^\infty({\mathbb R}^n)$, 
$C=C(n,\beta,\gamma)>0$. 
It is well known (see, e.g, \cite{HY}, \cite{HY2}, \cite{Me}) 
that this estimate with $\gamma=0$, 
together with the Duhamel principle, yields 
\smallskip

\noindent{\bf Lemma 6.1} {\it Let $n\geq 1$, $\delta>0$. 
Suppose that v solves the Cauchy problem 
$\Box v=G$ in $(0,T)\times{\mathbb R}^n$, 
with data $v(0)=f$, $\partial_t v(0)=g$. 
Then the estimate 
\begin{eqnarray}
& &
\|\langle r\rangle^{-(1/2)-\delta}\partial v\|
_{L^2((0,T)\times{\mathbb R}^n)}\\
& &
\leq
C
\bigl(
\|\nabla f\|_{L^2}+\|g\|_{L^2}
+
\|G\|_{L^1((0,T);L^2({\mathbb R}^n))}
\bigr)\nonumber
\end{eqnarray}
holds.
}

For $n=1$ or $n\geq 3$, 
the estimate (6.26) with $G\equiv 0$ is proved by the multiplier method, 
as mentioned on page 7 of \cite{HY}. 
On the other hand, the proof of the Smith-Sogge estimate (6.25) builds on 
Fourier analysis, 
and it is valid for $n=2$ as well. 
This is why we start with (6.25) for the proof of (6.26).

Using (4.4) with $|a|=\mu-1$ and (6.26) with $\delta=1/2$, 
and then proceeding as in (5.3)-(5.6), 
we obtain for $T<T_0$
\begin{eqnarray}
& &
\sum_{|a|\leq \mu-1}
\|
\langle r\rangle^{-1}
\partial\Gamma^a u
\|_{L^2((0,T)\times{\mathbb R}^2)}\\
& &
\leq
CE_\mu^{1/2}(u(0))
+
C\sum_{{|b|+|c|+|d|}\atop{\leq\mu-1}}
\int_0^T
\bigl(
\|\partial\Gamma^b u^i(t)
\partial\Gamma^c u^j(t)
\partial^2\Gamma^d u^l(t)\|_{L^2}\nonumber\\
& &
\hspace{5.0cm}
+\|\partial\Gamma^b u^i(t)
\partial\Gamma^c u^j(t)
\partial\Gamma^d u^l(t)\|_{L^2}
\bigr)dt\nonumber\\
& &
\leq
C
\bigl(
\varepsilon
+
\varepsilon^3
\int_0^T(1+t)^{-1}dt
\bigr)
\leq
C\varepsilon
\log (2+T),\nonumber
\end{eqnarray}
which is the estimate we have seeked for. 
\smallskip

\noindent{\it ${\dot H}^{1/2}$ estimate}. 
We use the following basic estimate. 
\smallskip

\noindent{\bf Lemma 6.2} {\it Let $v$ be the solution 
to the Cauchy problem 
$\Box v=G$ in $(0,T)\times{\mathbb R}^2$ 
with data $(v(0),\partial_t v(0))=(f,g)$. 
Then there holds that 
\begin{eqnarray}
& &
\||D_x|^{1/2}v(t)\|_{L^2({\mathbb R}^2)}
+
\||D_x|^{-1/2}\partial_t v(t)\|_{L^2({\mathbb R}^2)}\\
& &
\leq
C
\bigl(
\||D_x|^{1/2}f\|_{L^2}
+
\||D_x|^{-1/2}g\|_{L^2}
+
\|G\|_{L^1((0,t);L^{4/3}({\mathbb R}^2))}
\bigr).\nonumber
\end{eqnarray}
}

The proof is elementary, and we may omit it. 

Recall the definition (2.2) of $N_\mu(u(t))$. 
Using (4.4) with $|a|\leq\mu$, 
proceeding as in (5.3)-(5.6), 
and applying the H\"older inequality and (4.13) with $p=8$, 
we get for $t<T_0$
\begin{eqnarray}
& &
N_\mu(u(t))\\
& &
\leq
N_\mu(u(0))
+
C\int_0^t
\bigl(
\langle \tau\rangle^{-3/8}
\bigr)^2
\bigl(
M_\mu^{1/2}(u(\tau))E_\mu^{1/4}(u(\tau))
\bigr)^2
E_\kappa^{1/2}(u(\tau))d\tau\nonumber\\
& &
\leq
N_\mu(u(0))
+
C\varepsilon^2
\int_0^t
\langle \tau\rangle^{-3/4}
E_\kappa^{1/2}(u(\tau))d\tau.\nonumber
\end{eqnarray}

Now we are ready to complete the proof of our main theorem. 
Since we know $E_\mu^{1/2}(u(t))<2\varepsilon$ $(0\leq t<T_0)$ 
for a sufficiently small $\varepsilon$ such that 
$2\varepsilon\leq\varepsilon_0$ for 
$\varepsilon_0$ in (5.7), 
we get from (6.10)
\begin{equation}
{\tilde E}_\kappa(u(t))
\leq
{\tilde E}_\kappa(u(0))
\langle t\rangle^{C\varepsilon^2},
\end{equation}
which, combined with (6.29), yields
\begin{equation}
N_\mu(u(t))
\leq
N_\mu(u(0))
+
C\varepsilon^2
{\tilde E}_\kappa^{1/2}(u(0))
\langle t\rangle^{(1/4)+C\varepsilon^2}.
\end{equation}
Inserting (6.30) and (6.31) into (6.24) 
and using the obvious inequality 
${\tilde E}_\mu(u(t))\leq{\tilde E}_\kappa(u(t))$ as well, 
we have
\begin{eqnarray}
& &
{\tilde E}_\mu (u(t))\\
& &
\leq
{\tilde E}_\mu (u(0))
+
C\int_0^t
\bigl(
\langle\tau\rangle^{-(3/2)+C\varepsilon^2}
{\tilde E}_\kappa(u(0))\nonumber\\
& &
\hspace{3.4cm}
+
\langle\tau\rangle^{-1+C\varepsilon^2}
{\tilde E}_\kappa^{1/2}(u(0))
\sum_{|d|\leq\mu-1}
\|\langle r\rangle^{-1}
\partial\Gamma^d u(\tau)\|_{L^2({\mathbb R}^2)}\nonumber\\
& &
\hspace{3.4cm}
+\langle\tau\rangle^{-(3/2)+C\varepsilon^2}
{\tilde E}_\kappa^{1/2}(u(0))N_\mu(u(0))\nonumber\\
& &
\hspace{3.4cm}
+\langle\tau\rangle^{-(5/4)+C\varepsilon^2}
C\varepsilon^2{\tilde E}_\kappa(u(0))
\bigr)
{\tilde E}_\mu (u(\tau))d\tau\nonumber
\end{eqnarray}
for $0\leq t<T_0$. By the Gronwall inequality, we obtain 
\begin{equation}
{\tilde E}_\mu (u(t))
\leq
{\tilde E}_\mu (u(0))
\exp
\bigl(
C{\tilde E}_\kappa(u(0))
+
C{\tilde E}_\kappa^{1/2}(u(0))N_\mu(u(0))
\bigr)
\end{equation}
for $0\leq t<T_0$. 
Note that here we have used (6.27) together with 
the useful technique of 
dyadic decomposition of the interval $(0,T_0)$, 
as in page 363 of \cite{So}, page 408 of \cite{MNS}, 
page 13 of \cite{HY}, and page 11 of \cite{HY2}.

Recalling (6.5) and taking the size condition (3.6) into account, 
we see that there exists a constant $A$ such that 
\begin{eqnarray}
& &
E_\mu^{1/2}(u(t))\\
& &
\leq
2E_\mu^{1/2}(u(0))
\exp
\bigl(
AE_\kappa^{1/2}(u(0))
(E_\kappa^{1/2}(u(0))+N_\mu(u(0)))
\bigr)
<2\varepsilon\nonumber
\end{eqnarray}
for $0\leq t<T_0$. 
The last inequality proves that 
the norm $E_\mu^{1/2}(u(t))$ 
is strictly smaller than $2\varepsilon$ 
on the closed interval $[0,T_0]$. The proof 
of the main theorem has been completed.
$\hfill\square$

\bigskip

\noindent{\em Acknowledgement.} 
The author is grateful to Professor Kazuyoshi Yokoyama 
for a helpful discussion on the subject of this paper. 
He was supported in part by 
the Grant-in-Aid for Scientific Research (C) (No.\,23540198), 
Japan Society for the Promotion of Science (JSPS).

\begin{flushleft}
Kunio Hidano\\
Department of Mathematics\\
Faculty of Education\\
Mie University\\
1577 Kurima-machiya-cho, Tsu\\
Mie 514-8507 Japan\\
\end{flushleft}
\end{document}